\documentclass[leqno]{article}
\usepackage{amsmath}
\usepackage{amssymb}
\usepackage{dsfont}
\usepackage[usenames]{color}
\usepackage[colorlinks]{hyperref}
\usepackage{maple2e}
\def\ds{\displaystyle}
\begin{document}

\title{\bf Extremal orders of compositions of certain arithmetical functions}
\author{{\sc J\'ozsef S\'andor} and {\sc L\'aszl\'o T\'oth}}
\date{}
\maketitle

\centerline{Integers: Electronic Journal of Combinatorial Number
Theory, {\bf 8} (2008), \#A34} \vskip4mm

\begin{abstract} We study the exact extremal orders of compositions $f(g(n))$
of certain arithmetical functions, including the functions
$\sigma(n)$, $\phi(n)$, $\sigma^*(n)$ and $\phi^*(n)$, representing
the sum of divisors of $n$, Euler's function and their unitary
analogues, respectively. Our results complete, generalize and refine
known results.
\end{abstract}



\vskip2mm {\bf 1. Introduction}

\vskip1mm Let $\sigma(n)$, $\phi(n)$ and $\psi(n)$ denote -- as
usual -- the sum of divisors of $n$, Euler's function and the
Dedekind function, respectively, where $\psi(n)=n\prod_{p\mid n}
(1+1/p)$.

Extremal orders of the composite functions $\sigma(\phi(n))$, $\phi(\sigma(n))$, $\sigma(\sigma(n))$, $\phi(\phi(n))$, $\phi(\psi(n))$, $\psi(\phi(n))$, $\psi(\psi(n))$
% of these arithmetical functions
were investigated by {\sc L. Alaoglu} and {\sc P. Erd\H os}
\cite{AE1944}, {\sc A. M\c{a}kowski} and {\sc A. Schinzel}
\cite{MS1964}, {\sc J. S\'andor} \cite{S2005}, {\sc F. Luca} and {\sc
C. Pomerance} \cite{LP2002}, {\sc J.-M. de Koninck} and {\sc F. Luca}
\cite{KL2007}, and others.

For example, in paper \cite{MS1964} it is shown that
\begin{equation}
\liminf_{n\to \infty} \frac{\sigma(\sigma(n))}{n}=1, \label{sigsig}
\end{equation}
\begin{equation}
\limsup_{n\to \infty} \frac{\phi(\phi(n))}{n}=\frac1{2}, \label{phiphi}
\end{equation}
while in paper \cite{LP2002} the result
%\vskip1mm {\sc F. Luca}, {\sc C. Pomerance}
\begin{equation}
\limsup_{n\to \infty} \frac{\sigma(\phi(n))}{n\log \log n}
=e^{\gamma} \label{sigphi}
\end{equation}
is proved, where $\gamma$ is Euler's constant.

It is the aim of the present paper to extend the study of exact
extremal orders to other compositions $f(g(n))$ of arithmetical
functions, considering also the functions $\sigma^*(n)$ and
$\phi^*(n)$, representing the sum of unitary divisors of $n$ and the
unitary Euler function, respectively. Recall that $d$ is a unitary
divisor of $n$ if $d\mid n$ and $(d, n/d) = 1$. The functions
$\sigma^*(n)$ and $\phi^*(n)$ are multiplicative and if $n =
p_1^{a_1}\cdots p_r^{a_r}$ is the prime factorization of $n>1$, then
\begin{equation}
\sigma^*(n) =(p_1^{a_1}+ 1)\cdots (p_r^{a_r} + 1), \quad \phi^*(n) = (p_1^{a_1}-1)\cdots (p_r^{a_r}-1). \label{unitary}
\end{equation}

Note that $\sigma^*(n)=\sigma(n)$, $\phi^*(n)=\phi(n)$ for all
squarefree $n$ and that for every $n\ge 1$,
\begin{equation}
\phi(n)\le \phi^*(n)\le n \le \sigma^*(n)\le \psi(n)\le \sigma(n). \label{unitineq}
\end{equation}

We give some general results which can be applied easily also for
other special functions. Our results complete, generalize and refine
known results. They are stated in Section 2, their proofs are given
in Section 3. Some open problems are formulated in Section 4.

%\newpage
\vskip2mm

{\bf 2. Main results}

\vskip1mm {\bf Theorem 1.} {\sl  Let $f$ be an arithmetical
function. Assume that

(i) $f$ is integral valued and $f(n)\ge 1$ for every $n\ge 1$,

(ii) $f(n)\le n$ for every sufficiently large $n$ ($n\ge n_0$),

(iii) $f(p)=p-1$ for every sufficiently large prime $p$ ($p\ge
p_0$).

Then
\begin{equation}
\limsup_{n\to \infty} \frac{\sigma(f(n))}{n\log \log n}
=\limsup_{n\to \infty} \frac{\sigma(f(n))}{f(n)\log \log f(n)}=
e^{\gamma}, \label{sigf}
\end{equation}
\begin{equation}
\limsup_{n\to \infty} \frac{\psi(f(n))}{n\log \log n}= \limsup_{n\to
\infty} \frac{\psi(f(n))}{f(n)\log \log f(n)} =\frac{6}{\pi^2}
e^{\gamma}, \label{psif}
\end{equation}
\begin{equation}
\limsup_{n\to \infty} \frac{\sigma(f(n))}{\phi(f(n))(\log \log
n)^2}= \limsup_{n\to \infty} \frac{\sigma(f(n))}{\phi(f(n))(\log
\log f(n))^2} =e^{2\gamma}, \label{sigphif}
\end{equation}
\begin{equation}
\limsup_{n\to \infty} \frac{\psi(f(n))}{\phi(f(n))(\log \log n)^2}=
\limsup_{n\to \infty} \frac{\psi(f(n))}{\phi(f(n))(\log \log
f(n))^2} =\frac{6}{\pi^2} e^{2\gamma}. \label{psiphif}
\end{equation}}

\vskip1mm Theorem 1 can be applied for $f(n)=\phi(n)$ and
$f(n)=\phi^*(n)$, the unitary Euler function.
For example, \eqref{sigf} and \eqref{psif} give
\begin{equation}
\limsup_{n\to \infty} \frac{\sigma(\phi^*(n))}{n\log \log n}
=e^{\gamma}, \label{sigphi*}
\end{equation}
\begin{equation}
\limsup_{n\to \infty} \frac{\psi(\phi(n))}{n\log \log n}
=\frac{6}{\pi^2} e^{\gamma}. \label{psiphi}
\end{equation}

The weaker result $\ds \limsup_{n\to \infty}
\frac{\psi(\phi(n))}{n}=\infty$ is proved in \cite{S2005}.

\vskip1mm Figure 1 is a plot of the functions $\sigma(\phi^*(n))$
and $e^{\gamma}n\log \log n$ for $10\le n\le 10\, 000$.

\vskip1mm {\bf Theorem 2.} {\sl Let $g$ be an arithmetical function.
Assume that

(i) $g$ is integral valued and $g(n)\ge 1$ for every $n\ge 1$,

(ii) $g(n)\ge n$ for every sufficiently large $n$ ($n\ge n_0$),

(iii) either $g(p)=p+1$ for every sufficiently large prime $p$
($p\ge p_0$), or $g$ is multiplicative  and $g(p)=p$ for every
sufficiently large prime $p$ ($p\ge p_0$).

Then
\begin{equation}
\liminf_{n\to \infty} \frac{\phi(g(n))\log \log n}{n} =
\liminf_{n\to \infty} \frac{\phi(g(n))\log \log g(n)}{g(n)}=
e^{-\gamma}. \label{phig}
\end{equation}}

\vskip1mm Theorem 2 applies for $g(n)=\sigma(n), \sigma^*(n),
\psi(n), \sigma^{(e)}(n)$, where $\sigma^{(e)}(n)$ represents the
sum of exponential divisors of $n$. We have for example
\begin{equation}
\liminf_{n\to \infty} \frac{\phi(\sigma(n))\log \log n}{n} =e^{-\gamma}. \label{phisig}
\end{equation}

\vskip1mm Remark that according to a result of {\sc L. Alaoglu} and
{\sc P. Erd\H os} \cite{AE1944}, $\ds \lim_{n\to \infty}
\frac{\phi(\sigma(n))}{n} =0$ on a set of density $1$.

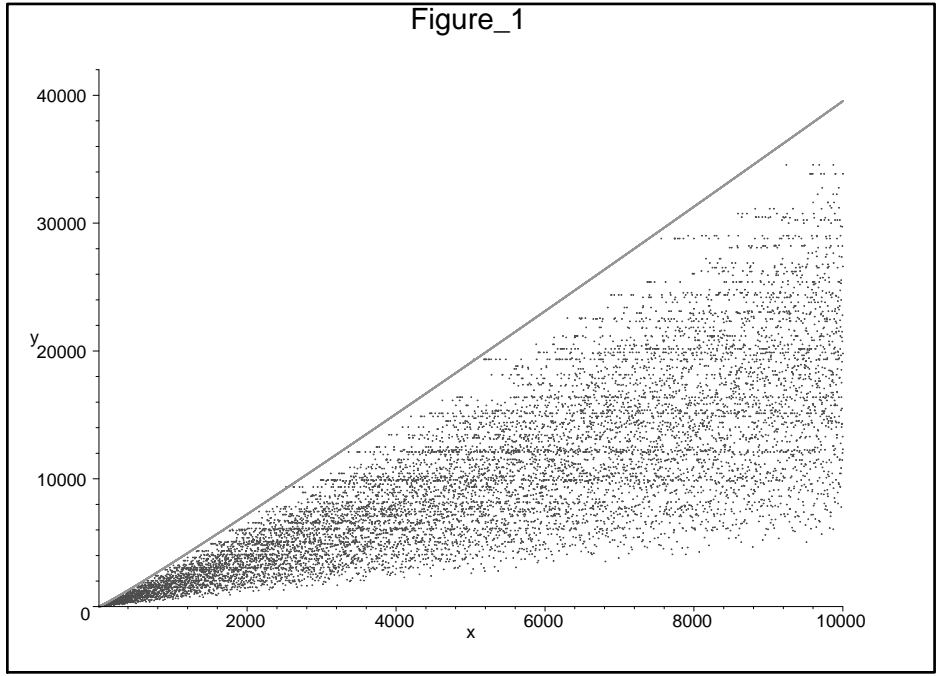
\begin{maplegroup}
\mapleresult
\begin{center}
\mapleplot{Fig101.eps} % Figure 1
\end{center}
\end{maplegroup}

\vskip1mm Theorems 1 and 2 can be generalized as follows. If
$f(n)\ge 1$ is an integer valued arithmetical function let $f_k(n)$
denote its $k$-fold iterate, i. e. $f_0(n)=n$, $f_1(n)=f(n)$, ...,
$f_k(n)=f(f_{k-1}(n))$.

\vskip1mm {\bf Theorem 3.} {\sl Let $f$ be an arithmetical function.
Suppose that

(i) $f$ is integral valued and $1\le f(n)\le n$ for every $n\ge 1$,

(ii) $f(p)=p-1$ for every prime $p$,

(iii) for every $s,t\ge 1$ if $s\mid t$, then $f(s)\mid f(t)$.

Then for every $k\ge 0$,
\begin{equation}
\limsup_{n\to \infty} \frac{\sigma(f_k(n))}{f_k(n)\log \log n}
=e^{\gamma}. \label{sigf_k/f_k}
\end{equation}}

\vskip1mm Theorem 3 applies for $f(n)=\phi(n)$, $f(n)=(p_1-1)\cdots
(p_r-1)$, $f(n)=(p_1-1)^{a_1}\cdots (p_r-1)^{a_r}$, where
$n=p_1^{a_1}\cdots p_r^{a_r}$.

\vskip1mm {\bf Theorem 4.} {\sl Let $g$ be an arithmetical function.
Suppose that

(i) $g$ is integral valued and $g(n)\ge n$ for every $n\ge 1$,

(ii) $g(p)=p+1$ for every prime $p$,

(iii) for every $s,t\ge 1$ if $s\mid t$, then $g(s)\mid g(t)$.

Then for every $k\ge 0$,
\begin{equation}
\liminf_{n\to \infty} \frac{\phi(g_k(n))\log \log n}{g_k(n)}
=e^{-\gamma}. \label{phig_k/g_k}
\end{equation}}

\vskip1mm Theorem 4 applies for $g(n)=\psi(n)$, $g(n)=(p_1+1)\cdots
(p_r+1)$, $g(n)=(p_1+1)^{a_1}\cdots (p_r+1)^{a_r}$, where
$n=p_1^{a_1}\cdots p_r^{a_r}$.

For $f(n)=\phi(n)$ and $g(n)=\psi(n)$ we have for every $k\ge 0$,
\begin{equation}
\limsup_{n\to \infty} \frac{\sigma(\phi_k(n))}{\phi_k(n)\log \log n}
=e^{\gamma}, \label{phi_k}
\end{equation}
\begin{equation}
\liminf_{n\to \infty} \frac{\phi(\psi_k(n))}{\psi_k(n)\log \log n}
=e^{-\gamma}. \label{psig_k}
\end{equation}

Compare Theorems 1--4 with the following deep results:

-- for $k\ge 2$ the normal order of $\ds
\frac{\sigma_k(n)}{\sigma_{k-1}(n)}$ is $ke^{\gamma}\log \log \log
n$, i.e. $\sigma_k(n)\sim ke^{\gamma} \sigma_{k-1}(n) \log \log \log
n$ on a set of density $1$, cf. {\sc P. Erd\H os} \cite{E1967},

-- for $k\ge 1$ the normal order of $\ds
\frac{\phi_k(n)}{\phi_{k+1}(n)}$ is $ke^{\gamma}\log \log \log n$,
proved by {\sc P. Erd\H os}, {\sc A. Granville}, {\sc C. Pomerance}
and {\sc C. Spiro} \cite{EGPS1990}.

-- the normal order of $\ds \frac{\phi(\sigma(n))}{\sigma(n)}$ is
$e^{-\gamma}/\log \log \log n$ and the normal order of $\ds
\frac{\sigma(\phi(n))}{\phi(n)}$ is $e^{\gamma}\log \log \log n$,
see {\sc L. Alaoglu} and {\sc P. Erd\H os} \cite{AE1944}.

\vskip1mm Note that the average orders of $\phi(n)/\phi_2(n)$ and
$\phi_2(n)/\phi(n)$ were investigated by {\sc R. Warlimont}
\cite{W2001}.

\vskip1mm {\bf Theorem 5.} {\sl Let $h(n)$ be an arithmetical
function such that $n\le h(n)\le \sigma(n)$ for every sufficiently
large $n$ ($n\ge n_0)$. Then
\begin{equation}
\liminf_{n\to \infty} \frac{h(\sigma(n))}{n} =1. \label{hsig}
\end{equation}}

\vskip1mm For $h(n)=\sigma(n)$ this is formula \eqref{sigsig}, for
$h(n)=\psi(n)$ it is due by {\sc J. S\'andor} \cite{S2005}, Theorem
3.30. Theorem 5 applies also for  $h(n)=\sigma^*(n),
\sigma^{(e)}(n)$.

\vskip1mm {\bf Theorem 6.} {\sl
\begin{equation}
\limsup_{n\to \infty} \frac{\phi(\phi^*(n))}{n}= \limsup_{n\to
\infty} \frac{\phi^*(\phi(n))}{n}= \limsup_{n\to
\infty} \frac{\phi^*(\phi^*(n))}{n}=1. \label{phiphi*}
\end{equation}}

Compare the results of \eqref{phiphi*} with \eqref{phiphi}.

\vskip1mm Figure 2 is a plot of the functions $\phi^*(\phi(n))$ and $n$ for
$1\le n\le 10\, 000$.

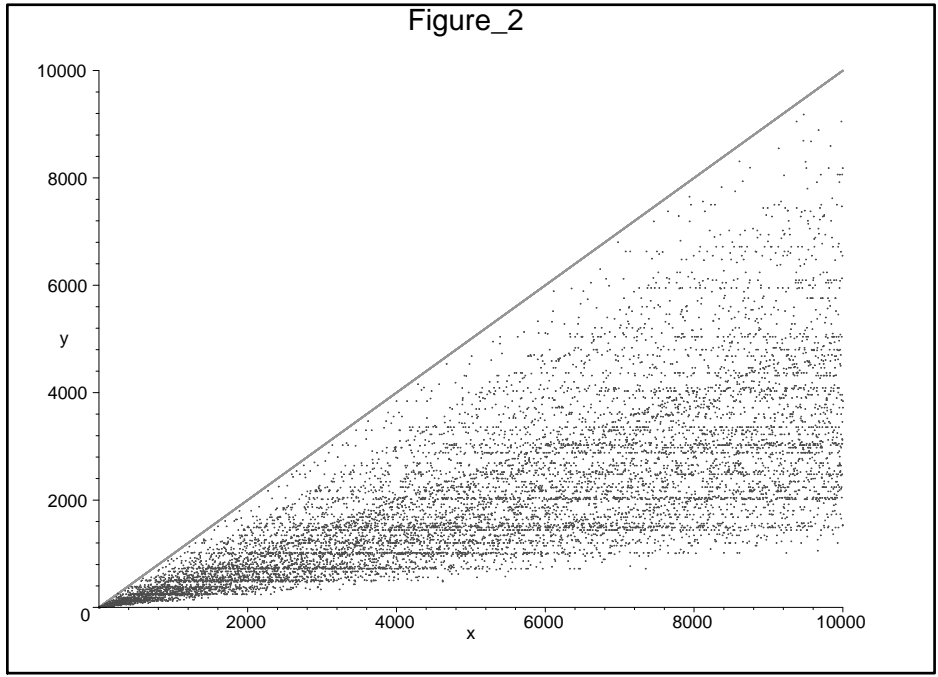
\begin{maplegroup}
\mapleresult
\begin{center}
\mapleplot{Fig2101.eps} %Figure 2
\end{center}
\end{maplegroup}

\vskip1mm Concerning $\phi^*(\phi^*(n))$ and $\sigma^*(\phi^*(n))$ we also prove:

\vskip1mm {\bf Theorem 7.}
\begin{equation}
\liminf_{n\to \infty} \frac{\phi^*(\phi^*(n))}{\log n \log \log n}>
0. \label{phi*phi*}
\end{equation}

\vskip1mm {\bf Theorem 8.} {\sl
\begin{equation}
\liminf_{n\to \infty} \frac{\sigma^*(\phi^*(n))}{n} \le \inf \left\{
\frac{\sigma^*(\phi^*(m/2))}{m/2} : 2\mid m, m\ne 2^{\ell}, \ell \ge
2 \right\}, \label{sig*phi*}
\end{equation}
\begin{equation}
\liminf_{n\to \infty} \frac{\sigma^*(\phi^*(n))}{n} \le \frac1{4} +
\varepsilon, \label{corol}
\end{equation}
where $\ds \varepsilon = \frac{3}{4(2^{32}-1)}\approx 0.17\cdot
10^{-9}$. }

%%%%%%%%%%%%%%%%%%%%%%%%%%%%%%%%%%%%%%%%%%%%%%%%%%%%%%%%%%%%%%%%%%%%%%%%%%%%%%%%%%%%%%%%%%%%%%%

\vskip2mm

{\bf 3. Proofs}

\vskip1mm The proofs of Theorems 1 and 2 are similar to the proof of
\eqref{sigphi} given in \cite{LP2002}, using a simple argument based
on Linnik's theorem, which states that if $(k,\ell)=1$, then there
exists a prime $p$ such that $p\equiv \ell$ (mod $k$) and $p\ll
k^{c}$, where $c$ is a constant (one can take $c\le 11$).

\vskip1mm {\bf Proof of Theorem 1.} To obtain the maximal orders of
the functions $\sigma(n)/n$, $\psi(n)/n$, $\sigma(n)/\phi(n)$ and
$\psi(n)/\phi(n)$, which are needed in the proof, we apply the
following result of {\sc L. T\'oth} and {\sc E. Wirsing}, see
\cite{TW2003}, Corollary 1:

If $F$ is a nonnegative real-valued multiplicative arithmetic
function such that for each prime $p$,

a) $\rho(p):=\sup_{\nu \ge 0} F(p^{\nu})\le (1-1/p)^{-1}$, and

b) there is an exponent $e_p=p^{o(1)}$ satisfying $F(p^{e_p})\ge
1+1/p$,

then
\[
\displaystyle \limsup_{n\to \infty} \frac{F(n)}{\log \log
n}=e^{\gamma}\prod_p \left(1-\frac1{p}\right) \rho(p).
\]

For $F(n)=\sigma(n)/n$ (with $\rho(p)=(1-1/p)^{-1}$, $e_p=1$),
$F(n)=\psi(n)/n$ (with $\rho(p)=1+1/p$, $e_p=1$),
$F(n)=\sqrt{\sigma(n)/\phi(n)}$ (with $\rho(p)=(1-1/p)^{-1}$,
$e_p=1$) and $F(n)=\sqrt{\psi(n)/\phi(n)}$ (with
$\rho(p)=\sqrt{(p+1)/(p-1)}$, $e_p=1$), respectively, we obtain
\begin{equation}
\limsup_{n\to \infty} \frac{\sigma(n)}{n \log \log n} = e^{\gamma},
\label{sigma}
\end{equation}
\begin{equation}
\limsup_{n\to \infty} \frac{\psi(n)}{n \log \log n} =\frac{6}{\pi^2}
e^{\gamma}, \label{psi}
\end{equation}
\begin{equation}
\limsup_{n\to \infty} \frac{\sigma(n)}{\phi(n) (\log \log n)^2} =
e^{2\gamma}, \label{fracsigphi}
\end{equation}
\begin{equation}
\limsup_{n\to \infty} \frac{\psi(n)}{\phi(n) (\log \log n)^2}
=\frac{6}{\pi^2} e^{2\gamma}. \label{psiphi}
\end{equation}

Note that \eqref{sigma} is the result of {\sc T. H. Gronwall}
\cite{G1913}, \eqref{psiphi} is due to {\sc S. Wigert} \cite{W1931}
and \eqref{fracsigphi} is better than $\limsup_{n\to \infty}
\sigma(n)/\phi(n)= \infty$ given in \cite{So1977}.

Proof of \eqref{sigf}. Using assumption (ii),
\[
\ell_f:= \limsup_{n\to \infty} \frac{\sigma(f(n))}{n\log \log n}\le
\ell_f':= \limsup_{n\to \infty} \frac{\sigma(f(n))}{f(n)\log \log
f(n)}\le \limsup_{m\to \infty} \frac{\sigma(m)}{m\log \log m}
=e^{\gamma},
\]
according to \eqref{sigma}. For every $n$, let $p_n$ be the least
prime such that $p_n\ge p_0$ and $p_n\equiv 1$ (mod $n$). Here
$n\mid p_n-1$ and by Linnik's theorem $p_n\ll n^{c}$, so $\log \log
p_n \sim \log \log n$. Hence, using condition (iii),
\[
\frac{\sigma(f(p_n))}{p_n\log \log p_n}=
\frac{\sigma(p_n-1)}{p_n\log \log p_n}\sim
\frac{\sigma(p_n-1)}{(p_n-1) \log \log n} \ge \frac{\sigma(n)}{n\log
\log n},
\]
applying that if $s\mid t$, then $\sigma(s)/s =\sum_{d\mid s} 1/d\le
\sum_{d\mid t} 1/d = \sigma(t)/t$. We obtain that $\ell_f \ge
e^{\gamma}$, therefore $e^{\gamma}\le \ell_f \le \ell_f'\le
e^{\gamma}$, that is $\ell_f=\ell_f' =e^{\gamma}$.

Proofs of \eqref{psif}, \eqref{sigphif}, \eqref{psiphif}. Analogous
to the  method of above taking into account \eqref{psi},
\eqref{fracsigphi}, \eqref{psiphi} and that $s\mid t$ implies
$\psi(s)/s \le \psi(t)/t$, $\sigma(s)/\phi(s) \le
\sigma(t)/\phi(t)$, $\psi(s)/\phi(s) \le \psi(t)/\phi(t)$. $\Box$

\vskip1mm {\bf Proof of Theorem 2.} Similar to the proof of Theorem
1. We use the result of {\sc E. Landau} \cite{L1909},
\begin{equation}
\liminf_{n\to \infty} \frac{\phi(n) \log \log n}{n} = e^{-\gamma}.
\label{phi}
\end{equation}

By condition (ii) and using that the function $(\log \log x)/x$ is
decreasing for $x\ge x_0$,
\[
\ell_g:= \liminf_{n\to \infty} \frac{\phi(g(n))\log \log n}{n} \ge
\ell_g':= \liminf_{n\to \infty} \frac{\phi(g(n))\log \log
g(n)}{g(n)} \ge \liminf_{m\to \infty} \frac{\phi(m)\log \log
m}{m}=e^{-\gamma},
\]
according to \eqref{phi}.

Assume that $g(p)=p+1$ for every $p\ge p_0$. For every $n$, let
$q_n$ be the least prime such that $q_n\ge p_0$ and $q_n\equiv -1$
(mod $n$). Here $n\mid q_n+1$ and by Linnik's theorem $\log \log q_n
\sim \log \log n$. Hence
\[
\frac{\phi(g(q_n))\log \log q_n}{q_n} = \frac{\phi(q_n+1)\log \log
q_n}{q_n}\sim \frac{\phi(q_n+1)\log \log n}{q_n+1}\le
\frac{\phi(n)\log \log n}{n},
\]
applying that if $s\mid t$, then $\phi(s)/s \ge \phi(t)/t$. We
obtain that $e^{-\gamma}\ge \ell_g$, therefore $e^{-\gamma}\le
\ell_g' \le \ell_g \le e^{-\gamma}$, that is $\ell_g=\ell_g'
=e^{-\gamma}$.

Now suppose that $g$ is multiplicative and $g(p)=p$ for every prime
$p\ge p_0$. As it is known, in \eqref{phi} the liminf is attained
for $n=n_k=p_1\cdots p_k$, the product of the first $k$ primes, when
$k\to \infty$.

Therefore, since $g(n_k)=g(p_1\cdots p_k)=g(p_1)\cdots g(p_k)=
p_1\cdots p_k=n_k$,
\[
\lim_{k\to \infty} \frac{\phi(g(n_k))\log \log n_k}{n_k} =
\lim_{k\to \infty} \frac{\phi(n_k)\log \log n_k}{n_k}= e^{-\gamma}.
\ \ \Box \]

\vskip1mm {\bf Proof of Theorem 3.} By condition (i),
$f_2(n)=f(f(n))\le f(n)\le n$ and $f_k(n)\le n$ for every $k\ge 0$.
Therefore,
\[
\ell_k:= \limsup_{n\to \infty} \frac{\sigma(f_k(n))}{f_k(n) \log
\log n}\le \limsup_{n\to \infty} \frac{\sigma(f_k(n))}{f_k(n)\log
\log f_k(n)}\le \ell_0:= \limsup_{m\to \infty}
\frac{\sigma(m)}{m\log \log m} =e^{\gamma},
\]
by \eqref{sigma}, for every $k\ge 0$.

By (iii), if $s\mid t$, then $f(s)\mid f(t)$, $f_2(s)\mid f_2(t)$
and $f_k(s)\mid f_k(t)$ for every $k\ge 0$. Now let $k\ge 1$. If
$p_n$ is the least
prime such that $p_n\equiv 1$ (mod $n$), cf. proof of Theorem 1, then $n\mid p_n-1$
and $f_{k-1}(n)\mid f_{k-1}(p_n-1)$. Therefore, applying also (ii),
\[
\frac{\sigma(f_k(p_n))}{f_k(p_n)\log \log p_n} \sim
\frac{\sigma(f_{k-1}(p_n-1))}{f_{k-1}(p_n-1) \log \log n} \ge
\frac{\sigma(f_{k-1}(n))}{f_{k-1}(n) \log \log n}=\ell_{k-1},
\]

Hence $\ell_k\ge \ell_{k-1}$, and it follows $\ell_{k}\ge
\ell_{k-1}\ge ... \ge \ell_0$, $\ell_0 \le \ell_k \le \ell_0$,
$\ell_k=\ell_0=e^{\gamma}$. $\Box$

\vskip1mm {\bf Proof of Theorem 4.} Similar to the proof of Theorem
3. By condition (i), $g_2(n)=g(g(n))\ge g(n)\ge n$ and $g_k(n)\ge n$
for every $k\ge 0$. Therefore,
\[
L_k:= \liminf_{n\to \infty} \frac{\phi(g_k(n))\log \log n}{g_k(n)}
\ge \liminf_{n\to \infty} \frac{\phi(g_k(n))\log \log
g_k(n)}{g_k(n)}\ge L_0:= \limsup_{m\to \infty} \frac{\phi(m)\log
\log m}{m} =e^{-\gamma},
\]
by \eqref{phi}, for every $k\ge 0$.

By (iii), if $s\mid t$, then $g(s)\mid g(t)$, $g_k(s)\mid g_k(t)$
for every $k\ge 0$. Now let $k\ge 1$. If $q_n$ is the least prime
such that $q_n\equiv -1$ (mod $n$), cf. proof of Theorem 2, then
$n\mid q_n+1$ and $g_{k-1}(n)\mid g_{k-1}(q_n+1)$. Therefore,
applying also (ii),
\[
\frac{\phi(g_k(q_n))\log \log q_n}{g_k(q_n)} \sim
\frac{\phi(g_{k-1}(q_n+1)) \log \log n}{g_{k-1}(q_n+1)} \le
\frac{\phi(g_{k-1}(n))\log \log n}{g_{k-1}(n)}=L_{k-1},
\]

Hence $L_k \le L_{k-1}$, and it follows $L_k\le L_{k-1}\le ... \le
L_0$, $L_0 \le L_k \le L_0$, $L_k=L_0=e^{-\gamma}$. $\Box$

\vskip1mm {\bf Proof of Theorem 5.} By $h(n)\ge n$ we have
$h(\sigma(n))\ge \sigma(n)\ge n$, $h(\sigma(n))/n\ge 1$ ($n\ge
n_0$).

We use that for a fixed integer $a>1$ and with $p$ prime, for
$N(a,p)=\frac{a^p-1}{a-1}$ and for an arithmetical function
satisfying $\phi(n)\le F(n)\le \sigma(n)$ ($n\ge n_0$) one has
\begin{equation}
\lim_{p\to \infty} \frac{F(N(a,p))}{N(a,p)}=1, \label{FNa_p}
\end{equation}
cf. for ex. {\sc D. Suryanarayana} \cite{Su1977}.
%(also for $F(n)=\sigma(n)$)

For $p,q$ primes, $\sigma(q^{p-1})=\frac{q^p-1}{q-1}=N(q,p)$. We obtain, using \eqref{FNa_p},
\[
\frac{h(\sigma(q^{p-1}))}{q^{p-1}} =\frac{h(N(q,p))}{N(q,p))}\cdot
\frac{q^p-1}{q^{p-1}(q-1)} \to \frac{q}{q-1}, \ \text{as}  \ p\to
\infty, \] where $\frac{q}{q-1}<1+\epsilon$ for each $\epsilon >0$
if $q\ge q(\epsilon)$. $\Box$

\vskip1mm {\bf Proof of Theorem 6.} We have $\phi(n)\le n$ and
$\phi^*(n)\le n$ for all $n\ge 1$, hence $\phi(\phi^*(n))\le
\phi^*(n)\le n$. Similarly, $\phi^*(\phi^*(n))\le n$.

Let $n=2^p$, $p$ prime, then $\phi^*(n)=2^p-1$ and
\[ \frac{\phi(\phi^*(n))}{n} = \frac{\phi(2^p-1)}{2^p} = \frac{\phi(2^p-1)}{2^p-1}\cdot
\frac{2^p-1}{2^p}\to 1, \ \ p\to \infty, \] using \eqref{FNa_p} for $a=2$ and
$F(n)=\phi(n)$.

Similarly the relation for $\phi^*(\phi^*(n))$, using \eqref{FNa_p} for
$F(n)=\phi^*(n)$.

For $\phi^*(\phi(n))$ this can not be applied  and we need a special
treatment.

Let $\displaystyle M=\prod_{p\leq x}p^{a_p},$ where $a_p=\left\{
                                             \begin{array}{cl}
                                               [2\log x], & \hbox{if }
p<x^{1/2}, \\
                                               4, & \hbox{if }
p\in[x^{1/2},x]
                                             \end{array} (p \hbox{ prime}).
                                           \right.
$

Let $q$ be the least prime of the form $q\equiv M+1$ (mod $M^2$). By
Linnik's theorem one has $q \ll M^{c}$, where $c$ satisfies $c\le
11$.

Now, put $n=q$. Then $\phi(n)=q-1= M(1+kM)=MN$ for some $k$. Thus
$(M,N)=1$, so $N$ is free of prime factors $\leq x$. Since
$\phi^*$ is multiplicative, $\displaystyle
\frac{\phi^*(\phi(n))}{n}=\frac{\phi^*(M)}{M}\cdot
\frac{\phi^*(N)}{N}\cdot \frac{MN}{1+MN}$. Here $\displaystyle
\frac{MN}{1+MN}\to 1$, as $n=q\to \infty$, so it is sufficient to
study $\displaystyle\frac{\phi^*(M)}{M}$ and
$\displaystyle\frac{\phi^*(N)}{N}$. Clearly,
$\displaystyle\frac{\phi^*(M)}{M}=\prod_{p\leq
x}\frac{p^{a_p}-1}{p^{a_p}}=\prod_{p\leq
x}\left(1-\frac{1}{p^{a_p}}\right)$. If $p<x^{1/2},$ then
$p^{a_p}\geq 2^{[2\log x]}>x$ for sufficiently large $x$. Otherwise,
$p^{a_p}\geq (x^{1/2})^4=x^2>x$ again. So $p^{a_p}>x$ anyway,
implying that
\begin{equation} \label{phi*M}
\frac{\varphi^*(M)}{M}>\left(1-\frac{1}{x}\right)^{\pi(x)}=1+\mathcal{O}\left(\frac{1}{\log
x}\right).
\end{equation}
%as $\displaystyle\pi(x)=\mathcal{O}\left(\frac{x}{\log x}\right).$
%($\pi(x)=$ number of primes $\leq x$)

Remark that $\displaystyle M<\prod_{p<x^{1/2}}p^{2\log x}\cdot
\prod_{p\leq x} p^4<\exp\left(\mathcal{O}(x^{1/2}\log
x+x)\right)=\exp\big(\mathcal{O}(x)\big)$ by the well-known fact:
$\displaystyle \prod_{p\leq a} p = e^{\mathcal{O}(a)}$. From $q\ll
M^{c'}$ and $M<\exp\big(\mathcal{O}(x)\big)$, by $N\ll M^{10}$ it
follows also that
\begin{equation}
N <\exp\big(\mathcal{O}(x)\big). \label{N}
\end{equation}

Let now $\displaystyle N=\prod_{i=1}^k q_i^{b_i}$ be the prime
factorization of $N$. We have $\displaystyle\log N=\sum_{i=1}^k
b_i\log q_i>(\log x)\sum_{i=1}^k b_i$, as $q_i>x$ for all
$1\leq i \leq k$. Here $\displaystyle\sum_{i=1}^k b_i\geq k,$ thus
$\displaystyle k<\frac{\log N}{\log x}\ll\frac{x}{\log x}$ by
\eqref{N}. Thus
\begin{equation} \label{phi*N}
\frac{\phi^*(N)}{N}=\prod_{i=1}^k\left(1-\frac{1}{q_i^{b_i}}\right)> \left(1-\frac{1}{x}\right)^k \geq
\left(1-\frac{1}{x}\right)^{\mathcal{O}(x/\log
x)}= 1+\mathcal{O}\left(\frac{1}{\log x}\right).
\end{equation}

By \eqref{phi*M} and \eqref{phi*N}, $\displaystyle
\frac{\phi^*(\phi(n))}{n}> 1+\mathcal{O}\left(\frac{1}{\log
x}\right)$ for sufficiently large $n.$ As
$n\ll \exp \big(\mathcal{O}(x)\big)$, we get $\log n\ll x$, so
$\displaystyle \frac{\phi^*(\phi(n))}{n}\to 1$, as $n=q\to
\infty.$

As $\displaystyle
\frac{\phi^*(\phi(n))}{n}\leq \frac{\phi(n)}{n}\leq 1$, the
proof is ready. $\Box$

\vskip1mm {\bf Proof of Theorem 7.} For all $n\ge 1$, $\phi^*(n)\ge
P(n)-1$, where $P(n)$ is the greatest prime factor of $n$. Let
$n=2^p$, $p$ prime, then $\phi^*(\phi^*(n))=\phi^*(2^p-1)\ge
P(2^p-1)-1$. Now we use the following result of {\sc P. Erd\H os}
and {\sc T. N. Shorey} \cite{ES1976}: $P(2^p -1)\ge c p \log p$  for
every prime $p$, where $c > 0$ is an absolute constant, and obtain
\begin{equation}
\frac{\phi^*(\phi^*(n))}{\log n \log \log n}\ge \frac{c p \log p- 1}
{p \log 2 (\log p+ \log \log 2)} \to \frac{c}{\log 2}, \quad p\to
\infty,
\end{equation}
and the result follows. $\Box$

\vskip1mm {\bf Proof of Theorem 8.} To prove \eqref{sig*phi*},
remark that if $2 \mid m$ and $m\ne 2^{\ell}$ ($\ell \ge 2$), then
$m/2$ is not a power of $2$, so $\phi^*(m/2)$ will be even (having
at least an odd prime divisor). Since $2\mid \phi^*(m/2)$, one can
write $\sigma^*(2\phi^*(m/2))<2 \sigma^*(\phi^*(m/2))$. Let $p$ be a
sufficiently large prime ($p>p_0$), then $(p,m/2)=1$ and obtain
\[
\frac{\sigma^*(\phi^*(mp/2))}{mp/2}=
\frac{\sigma^*((p-1)\phi^*(m/2))}{mp/2}\le \]
\[ \le \frac{\sigma^*((p-1)/2) \sigma^*(2\phi^*(m/2))}{mp/2}\le
\frac{\sigma^*((p-1)/2)}{p/2} \cdot
\frac{\sigma^*(\phi^*(m/2))}{m/2}
\]
by the above remark.

It is known that $\ds \frac{F((p-1)/2)}{(p-1)/2} \to 1$, as $p\to
\infty$, for $F(n)=\sigma(n)$, see \cite{MS1964} and it follows that
it holds also for $F(n)=\sigma^*(n)$ and obtain \eqref{sig*phi*}.

Now for \eqref{corol} let $m=4(2^{32}-1)=4F_0F_1F_2F_3F_4$ be $4$ times the
product of the known Fermat primes. Then $\phi^*(m/2)=
\phi^*(2F_0F_1F_2F_3F_4)=2^{1+2+4+8+16}=2^{31}$, $\ds
\frac{\sigma^*(\phi^*(m/2))}{m/2}=
\frac{2^{31}+1}{2(2^{32}-1)}=\frac1{4}+\varepsilon$, with the given
value of $\varepsilon$. $\Box$

%%%%%%%%%%%%%%%%%%%%%%%%%%%%%%%%%%%%%%%%%%%%%%%%%%%%%%%%%%%%%%%%%%%%%%%%%%%%%%%%%%%%%%%%%%%%%%%%%%%%%%%

\vskip1mm {\bf 4. Open problems}

\vskip1mm {\bf Problem 1.} Are the results of Theorem 1 valid if
$f(n)\le n$ for each $n\ge n_0$ and $f(p)=p$ for each prime $p\ge
p_0$?

Let $n=p_1^{\nu_1}\cdots p_r^{\nu_r} >1$ be an integer. An integer
$a$ is called regular (mod $n$) if there is an integer $x$ such that
$a^2x\equiv a$ (mod $n$). Let $\varrho(n)$ denote the number of
regular integers $a$ (mod $n$) such that $1\le a\le n$. Here
$\varrho(n)=(\phi(p_1^{\nu_1})+1)\cdots (\phi(p_r^{\nu_r})+1)$, in
particular $\varrho(p)=p$ for every prime $p$, cf. {\sc L. T\'oth}
\cite{T}.

Does Theorem 1 hold for  $f(n)=\varrho(n)$?

\vskip1mm {\bf Problem 2.} The method of proof of Theorems 1--4 does
not work in case of $\sigma^*(\phi(n))$ and $\sigma^*(\phi^*(n))$,
for example. We have
\[
\limsup_{n\to \infty} \frac{\sigma^*(\phi(n))}{n\log \log n}\le
\limsup_{n\to \infty} \frac{\sigma^*(\phi(n))}{\phi(n)\log \log
\phi(n)} \le \limsup_{n\to \infty} \frac{\sigma^*(n)}{n\log \log n}
=\frac{6}{\pi^2} e^{\gamma},
\]
cf. \cite{TW2003}, but the second part of the proof can not be
applied, because $n\mid m$ does not imply $\sigma^*(n)/n\le
\sigma^*(m)/m$.

What are the maximal orders $\sigma^*(\phi(n))$ and
$\sigma^*(\phi^*(n))$?

\vskip1mm Figure 3 is a plot of the function $\sigma^*(\phi(n))$ for
$1\le n\le 10\, 000$.

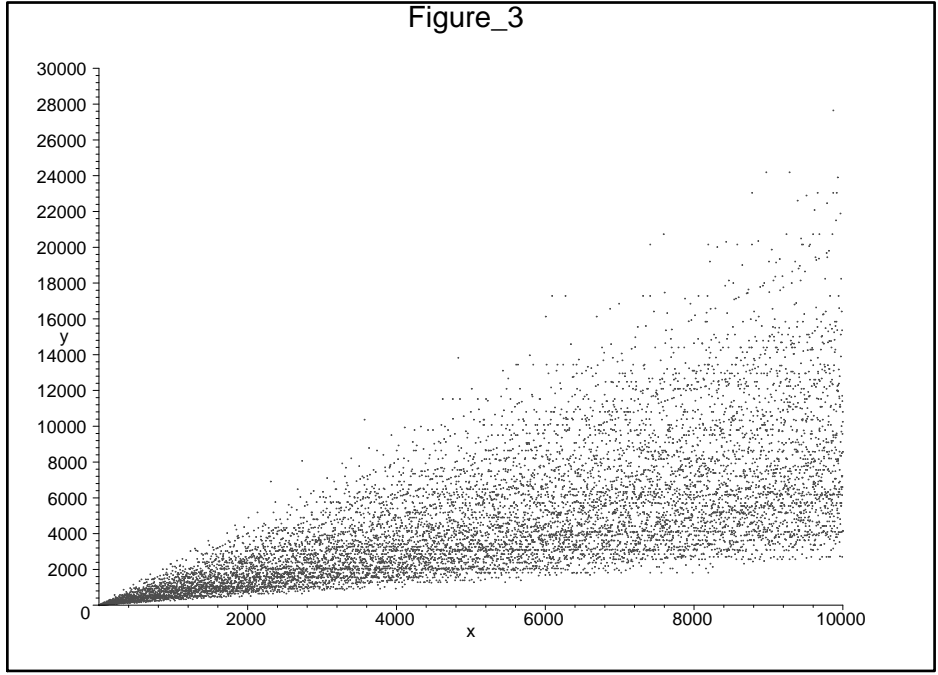
\begin{maplegroup}
\mapleresult
\begin{center}
\mapleplot{Fig301.eps} % Figure 3
\end{center}
\end{maplegroup}

\vskip1mm {\bf Problem 3.} Note that
\[
\limsup_{n\to \infty} \frac{\sigma^*(\sigma(n))}{n} = \limsup_{n\to
\infty} \frac{\sigma(\sigma^*(n))}{n} = \limsup_{n\to \infty}
\frac{\sigma^*(\sigma^*(n))}{n} =\infty,
\]
since for $n =n_k= p_1\cdots p_k$ (the product of the first $k$
primes),
\[ \frac{\sigma^*(\sigma(n_k))}{n_k}\ge
\frac{\sigma(n_k)}{n_k}=(1+1/p_1)\cdots (1+1/p_k)\to \infty, \quad
k\to \infty,
\]
similarly the other relations.

What are the maximal orders of $\sigma(\sigma^*(n))$,
$\sigma^*(\sigma(n))$,  $\sigma^*(\sigma^*(n))$?

\vskip1mm {\bf Problem 4.} Also,
\[
\liminf_{n\to \infty} \frac{\phi(\phi^*(n))}{n} = \liminf_{n\to
\infty} \frac{\phi^*(\phi(n))}{n} = \liminf_{n\to \infty}
\frac{\phi^*(\phi^*(n))}{n} =0,
\]
which follow at once by taking $n=n_k=p_1\cdots p_k$. Here
$\phi^*(\phi(n_k))= \phi^*((p_1-1)\cdots (p_k-1)) \le (p_1-1)\cdots
(p_k-1)-1$, hence
\[
\frac{\phi^*(\phi(n_k))}{n_k}\le \frac{(p_1-1)\cdots
(p_k-1)-1}{p_1\cdots p_k} < \left(1-\frac1{p_1}\right) \cdots
\left(1-\frac1{p_k}\right) \to 0, \ \ k\to \infty,
\]
similarly the other relations.

What are the minimal orders of $\phi(\phi^*(n))$, $\phi^*(\phi(n))$, $\phi^*(\phi^*(n))$?

%%%%%%%%%%%%%%%%%%%%%%%%%%%%%%%%%%%%%%%%%%%%%%%%%%%%%%%%%%%%%%%%%%%%%%%%%%%%%%%%%%%%%%%%%%%%%%%%%%%%%%%%

\vskip2mm {\bf 5. Maple notes}

The plots were produced using Maple. The functions $\sigma^*(n)$ and
$\phi^*(n)$ were generated by the following procedures:
\begin{verbatim}
sigmastar:= proc(n) local x, i: x:= 1: for i from 1 to nops(ifactors(n)[ 2 ]) do
p_i:=ifactors(n)[2][i][1]: a_i:=ifactors(n)[2][i][2];
x := x*(1+p_i^(a_i)): od: RETURN(x) end; # sum of unitary divisors

phistar:= proc(n) local x, i: x:= 1: for i from 1 to nops(ifactors(n)[ 2 ]) do
p_i:=ifactors(n)[2][i][1]: a_i:=ifactors(n)[2][i][2];
x := x*(p_i^(a_i)-1): od: RETURN(x) end; # unitary Euler function
\end{verbatim}

\vskip1mm {\bf Acknowledgements.} The authors wish to thank the
Referee for suggestions on improving earlier versions of Theorems 1
and 2, as well as for suggesting a correction for an initial version
of Theorems 3 and 4. The authors thank also Professor Florian Luca
for helpful correspondence.

\vskip3mm

{\sc J\'ozsef S\'andor}, Babe\c s-Bolyai University, Department of
Mathematics and Computer Sciences, Str. Kog\u alniceanu Nr. 1,
400084 Cluj-Napoca, Romania, E-mail: jjsandor@hotmail.com

\vskip2mm {\sc L\'aszl\'o T\'oth}, University of P\'ecs, Institute
of Mathematics and Informatics,  Ifj\'us\'ag u. 6, 7624 P\'ecs,
Hungary, E-mail: ltoth@ttk.pte.hu

\end{document}